\documentclass[graybox]{svmult}
\usepackage{mathptmx}
\usepackage{helvet}
\usepackage{courier}
\usepackage{makeidx}
\usepackage{graphicx}
\usepackage{multicol}
\usepackage[bottom]{footmisc}
\newcommand{\Teich}{\omega}

\newcommand{\Q}{{\mathbb{Q}}}

\newcommand{\Z}{{\mathbb{Z}}}

\usepackage{amsmath,amssymb,amsfonts}

\DeclareMathOperator{\Li}{Li}

\newcommand{\F}{\mathbb F}
\newcommand{\C}{\mathbb C}

\makeindex
\begin{document}
\title*{$P$-adic hypergeometrics}
\author{Fernando Rodriguez Villegas} 
\institute{Fernando Rodriguez Villegas \at Abdus Salam International Centre for
Theoretical Physics, Strada 
Costiera 11, 34151 Trieste, Italy,  
\email{villegas@ictp.it}
}
\maketitle
\abstract{We study classical hypergeometric series as  a $p$-adic function
  of its parameters inspired by a problem in the Monthly~\cite{Za}
  solved by D.~Zagier.} 
\bigskip

\noindent
\bigskip
\S1. The classical generalized hypergeometric series is defined by
\begin{equation}
\label{hyperg-series}
{}_rF_{r-1}\left[
 \begin{matrix}
 \alpha_1 & \ldots& \alpha_r \\
 \beta_1 & \ldots& \beta_{r-1}
\end{matrix}  \left.\right|\,t\right] =
\sum_{k\geq 0} \frac{(\alpha_1)_k\cdots (\alpha_r)_k}
{(\beta_1)_k\cdots (\beta_{r-1})_k}\frac{t^k}{k!}
\end{equation}
for $\alpha_j\in \C$ and $\beta_j\in\C \setminus\{0,-1,-2,\ldots\}$ .  If
$\alpha_r=-n$ for $n$ a non-negative integer the series 
terminates and we have
\begin{equation}
\label{hyperg-polyn}
{}_rF_{r-1}\left[
 \begin{matrix}
 \alpha_1 & \ldots& -n \\
 \beta_1 & \ldots& \beta_{r-1}
\end{matrix}  \left.\right|\,t\right] =
\sum_{k= 0}^n (-1)^k \binom n k \frac{(\alpha_1)_k\cdots (\alpha_{r-1})_k}
{(\beta_1)_k\cdots (\beta_{r-1})_k} t^k\;.
\end{equation}
One can show that for fixed $\alpha_i\in \Z_p,\beta_i\in
\Z_p\setminus\{0,-1,-2,\ldots\}$ and $|t|_p< 1$  
this yields  a convergent Mahler series  and hence 
a continuous function $f$ of the variable~$x:=n$ in $\Z_p$
\begin{equation}
\label{hyperg-polyn}
f(x):={}_rF_{r-1}\left[
 \begin{matrix}
 \alpha_1 & \ldots& -x \\
 \beta_1 & \ldots& \beta_{r-1}
\end{matrix}  \left.\right|\,t\right] =
\sum_{k\geq 0} (-1)^k \binom x k \frac{(\alpha_1)_k\cdots (\alpha_{r-1})_k}
{(\beta_1)_k\cdots (\beta_{r-1})_k} t^k\;.
\end{equation}
These functions seem very interesting and worthy of further
investigation. 

\bigskip
\noindent
\S2. It turns out that a special case of these functions appears in the
solution of an interesting Monthly problem~\cite{Za} solved by
D.~Zagier. The problem is to prove that
\begin{equation}
  \label{monthly-pblm}
  v_3\left(\sum_{k=0}^{n-1}\binom{2k}k\right)=  v_3\left(n^2\binom{2n}n\right),
\end{equation}
where $v_p$ denotes the $p$-adic valuation. Zagier does this by showing
that there is a continuous function $f_1: \Z_3 \longrightarrow -1+3\Z_3$
which interpolates the values
\begin{equation}
  \label{f-integer-values}
  f_1(n)=\frac{\sum_{k=0}^{n-1}\binom{2k}k}{n^2\binom{2n}n}, \qquad n=1,2,\ldots. 
\end{equation}
Considering the expansion
$$
f_1(n)=A + Bn+Cn^2+\cdots
$$
he goes further and conjectures, based on numerical evidence, that
$B=0$; moreover, he mentions
\begin{quote}
{\em
  Another interesting problem would be to evaluate in closed form the
  $3$-adic number $A$.}
\end{quote}

We prove that in fact 
\begin{equation}
\label{A-zeta}
A=-\frac32\zeta_3(2)
=2 + 3 + 2\cdot3^2 + 2\cdot3^6 + 3^7 + 2\cdot3^8 + 2\cdot3^9 + O(3^{10}),
\end{equation}
where $\zeta_3(s)$ is the Kubota-Leopoldt $3$-adic zeta function.

\bigskip
\noindent
\S3. The connection with zeta values is perhaps to be expected: in
general the Taylor coefficients of the functions of~\S1 involve
multiple polylogarithms.  In the specific case in question we have
  \begin{equation}
    \label{mahler-exp}
    f(n)=\sum_{k=1}^n\frac1{\binom{2k}k}\binom nk(-3)^{k-1}, \qquad
    f(n)=nf_1(n). 
  \end{equation}
If we expand in general
$$
f(x):= \sum_{k\geq1}\frac1{\binom{2k}k}\binom xkt^{k-1}=\sum_{n\geq
   0}b_n(t)x^n,
$$
then
\begin{equation}
b_n(t)=\frac1{(t+4)} \sum_{0\leq j_1 < j_2 < \cdots <
 j_n} 
\frac{(\frac t{t+4})^{j_n}}{(j_1+\tfrac12)(j_2+\tfrac12)\cdots(j_n+\tfrac12)}. 
\end{equation}

These multiple polylogarithms can be expressed in terms of usual
polylogarithms for small $n$. Trivially $b_0=0$. For $n=1$ we have
the following identity of power series in  $z=1-w$
  \begin{equation}
    \label{b_1-fmla}
b_1\left((w-w^{-1})^2\right)=
(w^2-w^{-2})^{-1}\log \left(w^2\right),
 \qquad  w=1-z.
\end{equation}
By plugging in a primitive third root of unity $\zeta_3 \in\C_3$ for
$w$ it follows that $3$-adically we have $b_1(-3)=0$. This shows that
in this case $f(x)$ is divisible by $x$ and we may consider
$f_1(x):=f(x)/x$~(see~\cite{Za}).

With some effort one can prove that as power series in $z$, with $w=1-z$, we have
\begin{equation}
    \label{b_2-fmla}
\begin{split}
b_2\left((w-w^{-1})^2\right) =(w^2-w^{-2})^{-1}
	[\Li_2(1-w^2)-\tfrac12 \Li_2(1-w^4)\\
	-\Li_2(1-w^{-2})+\tfrac12 \Li_2(1-w^{-4})],
\end{split}
\end{equation}
where $\Li_2$ is the standard dilogarithm function. 

Plugging in $w=\zeta_3\in \C_3$
into~\eqref{b_2-fmla} and using a result of Coleman~\cite{Co} we
obtain~\eqref{A-zeta}. The identity is the special case $p=3,r=1$ of the
following. Given a prime $p>2$ fix $\zeta_p\in \C_p$ a primitive $p$-th
root of unity.
\begin{theorem}
\label{L-values-thm}
i) The following limit exists
\begin{equation}
\label{A-lim}
A(\zeta_p):=\lim_{s\to \infty}
\frac{1}{\binom{2
    p^s}{p^s}p^{2s}}\sum_{k=0}^{p^s-1}\binom{2k}k(\zeta_p+\zeta_p^{-1})^{2(p^s-1-k)} 
\end{equation}

ii) Let $\Teich: \F_p^\times \rightarrow \C_p^\times$ be the Teichm\"uller
character. For $0<r<p-1$ we have 
\begin{equation}
\frac{1}{(\Teich^r(4)-2\Teich^r(2))}\sum_{i=1}^{p-1} \Teich(i)^{-r}
(\zeta_p^{2i}-\zeta_p^{-2i}) 
\, A(\zeta_p^i)=L_p(2,\Teich^{r-1}),
\end{equation}
where  $L_p$ is Kubota-Leopoldt's $p$-adic $L$-function.
\end{theorem}

We note in passing that 
$$
\lim_{s\rightarrow \infty}\binom {2 p^s}{p^s} = 2\prod_{k\geq 1}
\frac{\Gamma_p(2p^k)}{\Gamma_p(p^k)^2} 
$$
(see~\cite[\S6.3.4, ex. 16]{FRV}), where $\Gamma_p$ denotes the
$p$-adic gamma function.

The beauty of the expressions~\eqref{b_1-fmla} and~\eqref{b_2-fmla} is
that though their proof were obtained working over the complex numbers
they are identities of power series with rational coefficients and
hence also hold $p$-adically in an appropriate domain. Fortunately,
this domain includes the point were need to evaluate  for Zagier's
questions ($w=\zeta_3\in\C_3$).

For $n=3$ there is an expression for $b_3(t)$ in terms of
polylogarithms valid over the complex numbers, which is much more
difficult to obtain. For $n>3$ we do not expect $b_n(t)$ to reduce to
polylogarithms.

However, to apply this expression for $b_3$ to our $p$-adic setting
requires some form of analytic continuation. This we will achieve by 
delicate manipulations using Coleman's integration but the details
have not yet been fully carried out.

The expectation nevertheless is that for $p=3$ we should have that
$b_3(-3)$ is a simple multiple of $L_3(3,\chi_{-3})$. But
$L_3(s,\chi_{-3})$ is identically zero since $\chi_{-3}$ is odd! Hence
the constant $B$ of Zagier should vanish because it is a special value
of an $L$-function which happens to be identically zero.

\bigskip
\noindent
\S4.
We tested numerically to see if there are any other relations for
$b_n(t)$ and $p$-adic $L$-values and found only the following likely
identities: 

\begin{equation}
\Q_3:
\begin{cases}
b_4(-3) = -\frac{27}8\zeta_3(4)\\
b_6(-3) = -\frac{297}{32}\zeta_3(6)\\
\end{cases}
\Q_5:
\begin{cases}
b_2(-5) = 0\\
b_3(-5) = -\frac{25}{12}\zeta_5(3)\\
\end{cases}
\end{equation}
but we did not attempt to prove these. We pointed out above that $b_4(t)$
and $b_6(t)$ are not expected to be expressible in terms of
polylogarithms. Hence the connection of the observed identities for
$b_4(-3)$ and $b_6(-3)$ in $\Q_3$ appear to be less obvious than the
others.

\end{document}